%% file: main.tex
\RequirePackage{fix-cm}
\documentclass[smallcondensed]{svjour3}     
\smartqed  

\usepackage{graphicx}
\usepackage{amsmath,amssymb}
\usepackage{tikz-cd}
\usepackage[retainorgcmds]{IEEEtrantools}

\smartqed

\input{sc}

\begin{document}

\title{Homotopy types of gauge groups of $\PU(p)$-bundles over spheres
}


\author{Simon Rea
}


\institute{Department of Mathematical Sciences\\
University of Southampton, Highfield Campus\\
SO17 1BJ,  United Kingdom\\
              \email{S.Rea@soton.ac.uk}
}

\date{Received: date / Accepted: date}

\maketitle

\begin{abstract}
We examine the relation between the gauge groups of $\SU(n)$- and $\PU(n)$-bundles over $S^{2i}$, with $2\leq i\leq n$, particularly when $n$ is a prime. As special cases, for $\PU(5)$-bundles over $S^4$, we show that there is a rational or $p$-local equivalence $\cG_{2,k}\simeq_{(p)}\cG_{2,l}$ for any prime $p$ if, and only if, $(120,k)=(120,l)$, while for $\PU(3)$-bundles over $S^6$ there is an integral equivalence $\cG_{3,k}\simeq\cG_{3,l}$ if, and only if, $(120,k)=(120,l)$.
\end{abstract}

\section{Introduction}
\label{intro}

Let $G$ be a topological group and $X$ a space. The \emph{gauge group} $\cG(P)$ of a principal $G$-bundle $P$ over $X$ is defined as the group of $G$-equivariant bundle automorphisms of $P$ which cover the identity on $X$. A detailed introduction to gauge groups can be found in \cite{husemoller,piccinini1998conjugacy}.

The following problem is of interest: having fixed a topological group $G$ and a space $X$, classify the possible homotopy types of the gauge groups $\cG(P)$ of principal $G$-bundles $P$ over $X$.

Crabb and Sutherland showed \cite[Theorem 1.1]{crabb-sutherland} that if $G$ is a compact, connected, Lie group and $X$ is a connected, finite CW-complex, then the number of distinct homotopy types of $\cG(P)$, as $P$ ranges over all principal $G$-bundles over $X$, is finite. This is often in contrast with the fact that the number of isomorphisms classes of principal $G$-bundles over $X$ may be infinite. However, their methods did not lead to an enumeration of the classes of gauge groups. Classification results require a different kind of analysis.

The first classification result was obtained by Kono \cite{kono91} in 1991. Using the fact that principal $\SU(2)$-bundles over $S^4$ are classified by $k\in\Z\cong\pi_3(\SU(2))$ and denoting by $\cG_k$ the gauge group of the principal $\SU(2)$-bundle $P_k\to S^4$ corresponding to the integer $k$, Kono showed that there is a homotopy equivalence $\cG_k\simeq\cG_l$ if, and only if, $(12,k)=(12,l)$, where $(m,n)$ denotes the greatest common divisor of $m$ and $n$. It thus follows that there are precisely six homotopy types of $\SU(2)$-gauge groups over $S^4$.

In this paper, we examine how the close relationship between the groups $\SU(n)$ and $\PU(n)$ is reflected in the homotopy properties of the gauge groups of the corresponding bundles, particularly when $n$ is a prime. We do this by generalising certain results relating the classification of $\PU(n)$-gauge groups to that of $\SU(n)$-gauge groups from the paper \cite{so3} for the case $n=2$ (observe that $\PU(2)\cong\SO(3)$), and from \cite{hasui16} for the case $n=3$. 

Our first main result compares certain Samelson products on $\SU(p)$ and $\PU(p)$, with $p\geq 3$ a prime. As will be illustrated in Section \ref{sec:gauge}, the finiteness of the orders of Samelson products (in the appropriate groups of homotopy classes of maps) plays a crucial role in the homotopy classification of gauge groups. In Section \ref{sec:samels}, we will show the following.

\bt
\label{thm:pup-samels}
Let $p$ be an odd prime and let $2\leq i\leq p$. Let $\epsilon_i$ and $\delta_i$ denote generators of $\pi_{2i-1}(\PU(p))$ and $\pi_{2i-1}(\PU(p))$, respectively. The orders of the Samelson products $\langle \epsilon_i,1\rangle\cl S^{2i-1}\land \PU(p)\to\PU(p)$ and $\langle \delta_i,1\rangle\cl S^{2i-1}\land \SU(p)\to\SU(p)$, where $1$ denotes the appropriate identity map, coincide.
\et

Theorem \ref{thm:pup-samels} is the key ingredient for ``if'' direction of our classification results. The converse direction, on the other hand, will require that suitable homotopy invariants of the gauge groups be identified. In Section \ref{sec:gauge}, we introduce the notation $\cG_{i,k}(\PU(n))$, with $k\in\Z$, for gauge groups of $\PU(n)$-bundles over $S^{2i}$ in analogy with the notation used by Kono \cite{kono91} and others. In Section \ref{sec:hinv}, we give a sufficient condition for certain homotopy invariants of $\SU(n)$- and $\PU(n)$-gauge groups to coincide.

Our methods then allow us to deduce classification results for $\PU(p)$-gauge groups from the corresponding classification results for $\SU(p)$-gauge groups. As examples of applications of our results, we obtain the following complete classifications.

\bt
\label{thm:pu5}
For $\PU(5)$-bundles over $S^4$, it is the case that

(a) if $\cG_{2,k}(\PU(5))\simeq \cG_{2,l}(\PU(5))$, then $(120,k)=(120,l)$;

(b) if $(120,k)=(120,l)$, then $\cG_{2,k}(\PU(5))\simeq \cG_{2,l}(\PU(5))$ when localised rationally or at any prime. 
\et

\bt
\label{thm:pu3}
For $\PU(3)$-bundles over $S^6$, we have $\cG_{3,k}(\PU(3))\simeq\cG_{3,l}(\PU(3))$ if, and only if, $(120,k)=(120,l)$.
\et

We should note that in \cite{hasui16}, the $\PU(3)$-gauge group $\cG_{2,k}$ is shown to be homotopy equivalent to $\widehat{\cG}_{2,k}\times S^1$, where $\widehat{\cG}_{2,k}$ is a space whose homotopy groups are all finite. This allows the authors to apply Lemma \ref{lem:hamanaka} and obtain a classification result for $\cG_k$ that holds integrally. We expect the same result to apply more generally to gauge groups of $\PU(n)$-bundles over $S^{2n-2}$. However, there are currently no other cases, beside that of \cite{hasui16}, in which such a result would be applicable.

Finally, it is worth noting that, should any further classifications of gauge groups of $\SU(p)$-bundles over even-dimensional spheres be obtained, our results would provide the corresponding classification results for $\PU(p)$-gauge groups as immediate corollaries, provided the original results were arrived at via the standard methods. 

\section{Homotopy types of $\PU(n)$-gauge groups}
\label{sec:gauge}

The \emph{projective unitary group} $\PU(n)$ is defined as the quotient of $\U(n)$ by its centre
\bse
Z(\U(n))=\{\lambda I_n \mid \lambda\in \C \text{ and }|\lambda|=1\}\cong \U(1).
\ese
In marked contrast with the orthogonal case, there is a homotopy equivalence between $\PU(n)$ and $\PSU(n)$ for every $n$. For the purposes of this paper, one can therefore equivalently define $\PU(n)$ as the quotient of $\SU(n)$ by its centre
\bse
Z(\SU(n))=\{\lambda I_n \mid \lambda\in \C \text{ and }\lambda^n=1\}\cong \Z/n\Z.
\ese
We let $q\cl \SU(n)\to\PU(n)$ denote the quotient map corresponding to the latter definition, which we shall use throughout. Note that $q$ is an $n$-fold covering map and $\SU(n)$, being simply-connected, is the universal cover of $\PU(n)$.

\subsection{Notation for $\PU(n)$-gauge groups}

Let $n$ and $i$ be fixed integers such that $2\leq i\leq n$. As stated in the introduction, we are interested in the problem of classifying the homotopy types of the gauge groups $\cG(P)$ as $P$ ranges over all $\PU(n)$-bundles over $S^{2i}$.

If there exists a $\PU(n)$-equivariant bundle isomorphism $P\cong P'$, then conjugation by such an isomorphism yields a homeomorphism $\cG(P)\cong\cG(P')$. It therefore suffices to let $P$ range over a set of representatives of all the isomorphism classes of $\PU(n)$-bundles over $S^{2i}$.

Since $S^{2i}$ is paracompact, the set of isomorphism classes of principal $\PU(n)$-bundles over $S^{2i}$ is in bijection with the set $[S^{2i},\B\PU(n)]_\mathrm{free}$ of free homotopy classes of maps from $S^{2i}$ to $\B\PU(n)$, the classifying space of $\PU(n)$. As $\PU(n)$ is connected, $\B\PU(n)$ is simply-connected and hence there is a bijection 
\bse
[S^{2i},\B\PU(n)]_\mathrm{free}\cong[S^{2i},\B\PU(n)],
\ese
the right-hand side denoting the set of pointed homotopy classes of maps from $S^{2i}$ to $\B\PU(n)$. Furthermore, we have bijections
\begin{equation*}
[S^{2i},\B\PU(n)]\cong \pi_{2i}(\B\PU(n))\cong\pi_{2i-1}(\PU(n)) \cong \Z.
\end{equation*}
We can therefore introduce the following labelling for the gauge groups of $\PU(n)$-bundles over $S^{2i}$. Let $\epsilon_i\cl S^{2i-1}\to\PU(n)$ denote a generator of $\pi_{2i-1}(\PU(n))$. Each isomorphism class of $\PU(n)$-bundles is represented by the bundle $P_k\to S^{2i}$ induced by pulling back the universal $\PU(n)$-bundle along the classifying map $k\overline{\epsilon}_i\cl S^{2i}\to\B\PU(n)$, where $\overline{\epsilon}_i$ denotes the adjoint of $\epsilon_i$ and generates $\pi_{2i}(\B\PU(n))$.  
We let $\cG_{i,k}(\PU(n))$ (or simply $\cG_{i,k}$ when the context is clear) denote the gauge group of $P_k\to S^{2i}$. Mutatis mutandis, the notation $\cG_{i,k}(\SU(n))$ will also be used.

\subsection{Classification of the homotopy types of $\cG_{i,k}$}

By \cite{atiyah-bott} or \cite{gottlieb}, there is a homotopy equivalence $\B \cG_{i,k}\simeq\Map_{k}(S^{2i},\B\PU(n))$, the latter space being the $k$-th component of $\Map(S^{2i},\B\PU(n))$, meaning the connected component containing the classifying map $k\overline{\epsilon}_i$.

There is an evaluation fibration
\bse
\Map^*_k(S^{2i},\B\PU(n))\longrightarrow \Map_k(S^{2i},\B\PU(n))\xrightarrow{\ \ev\ } \B\PU(n),
\ese
where $\ev$ evaluates a map at the basepoint of $S^{2i}$ and the fibre is the $k$-th component of the pointed mapping space $\Map^*_k(S^{2i},\B\PU(n))$. This fibration extends to a homotopy fibration sequence
\bse
\cG_{i,k}\longrightarrow \PU(n)\longrightarrow\Map^*_k(S^{2i},\B\PU(n))\longrightarrow \B\cG_{i,k}\longrightarrow\B\PU(n),
\ese
where we used the equivalences $\B \cG_{i,k}\simeq\Map_{k}(S^{2i},\B\PU(n))$, $\Omega\B\PU(n)\simeq\PU(n)$, and $\Omega\B\cG_{i,k}\simeq\cG_{i,k}$.
Furthermore, by \cite{sutherland92} there is, for each $k\in\Z$, a homotopy equivalence 
\bse
\Map^*_k(S^{2i},\B\PU(n))\simeq \Map^*_0(S^{2i},\B\PU(n)).
\ese
The space on the right-hand side is denoted $\Omega_0^{2i}\B\PU(n)$ and is homotopy equivalent to $\Omega_0^{2i-1}\PU(n)$. We therefore have the following homotopy fibration sequence
\bse
\cG_{i,k}\longrightarrow \PU(n)\xrightarrow{\ \p_{i,k}\ }\Omega_0^{2i-1}\PU(n)\longrightarrow \B\cG_{i,k}\longrightarrow\B\PU(n),
\ese
which exhibits the gauge group $\cG_{i,k}$ as the homotopy fibre of the map $\p_{i,k}$. This is a key observation, as it suggests that the homotopy theory of the gauge groups $\cG_{i,k}$ depends on the maps $\p_{i,k}$. In fact, more is true.


By \cite[Theorem 2.6]{lang-ehp}, the adjoint of $\p_{i,k}\cl \PU(n)\to \Omega_0^{2i-1}\PU(n)$ is homotopic to the Samelson product $\langle k\epsilon_i,1\rangle\cl S^{2i-1}\land\PU(n)\to\PU(n)$, where 1 denotes the identity map on $\PU(n)$. As the Samelson product is bilinear, $\langle k\epsilon_i,1\rangle \simeq  k\langle \epsilon_i,1\rangle$, and hence, taking adjoints once more, $\p_{i,k}\simeq k \p_{i,1}$.

Thus, every one of the gauge groups is the homotopy fibre of the map $\p_{i,1}$ post-composed with the appropriate power map on $\Omega^{2i-1}_0\PU(n)$. If $\p_{i,1}$ can be determined to have finite order in $[\PU(n),\Omega_0^{2i-1}\PU(n)]$, this will have strong implications for the homotopy types of the $\cG_{i,k}$'s, as the following lemmas show.

\bl[Theriault {\cite[Lemma 3.1]{theriault-sp2}}]
\label{lem:theriault}
Let $X$ be a connected CW-complex and let $Y$ be an H-space with a homotopy inverse. Suppose that $f\in[X,Y]$ has finite order $m$. Then, for any integers $k,l\in\Z$ such that $(m,k)=(m,l)$, the homotopy fibres of $k f$ and $lf$ are homotopy equivalent when localised rationally or at any prime. \qed
\el

If, additionally, the homotopy groups of $Y$ are all finite, as in the case of $\PU(n)$-bundles over $S^{2n}$, where $Y= \Omega^{2n-1}_0\PU(n)$, the following stronger lemma applies.

\bl[Hamanaka, Kono {\cite[Lemma 3.2]{hamanaka-kono}}]
\label{lem:hamanaka}
Let $X$ be a connected CW-complex and let $Y$ be an H-space such that $\pi_j(Y)$ is finite for all $j$. Let $f\in[X,Y]$ be such that $mf\simeq*$ for some finite $m$ and let $k,l\in\Z$ satisfy $(m,k)=(m,l)$. Then, there exists a homotopy equivalence $h\cl Y\to Y$ satisfying $hkf\simeq lf$. \qed
\el

Note that the order of $\p_{i,1}$ coincides with the order of $\langle \epsilon_i,1\rangle$.

\section{Samelson products on $\PU(p)$}
\label{sec:samels}

Having fixed $n\geq 3$ and $2\leq i\leq n$, let $\delta_i\cl S^{2i-1}\to \SU(n)$ denote the generator of $\pi_{2i-1}(\SU(n))\cong \Z$ corresponding to the generator $\epsilon_i$ of $\pi_{2i-1}(\PU(n))$. That is, such that $q_*(\delta_i)=\epsilon_i$, where $q$ denotes the quotient map $q\cl \SU(n)\to \PU(n)$ introduced in Section \ref{sec:gauge}.

In this section, we wish to compare the orders of the Samelson products $\langle \delta_i,1\rangle$ and $\langle \epsilon_i,1\rangle$ on $\SU(n)$ and $\PU(n)$, respectively. First, observe that there is a commutative diagram
\bse
\label{diag}
\tag{$\star$}
\begin{tikzcd}[row sep=large,column sep=large]
S^{2i-1}\land \SU(n) \ar[r,"{\langle\delta_i,1\rangle}"] \ar[d,"1\land q"]& \SU(n)\ar[d,"q"]\\
S^{2i-1}\land \PU(n) \ar[r,"{\langle\epsilon_i,1\rangle}"] & \PU(n)
\end{tikzcd}
\ese
and recall the following property of the quotient map $q$.
\bl
\label{lem:p-equiv}
The quotient map $q\cl \SU(n)\to \PU(n)$ induces a $p$-local homotopy equivalence $\SU(n)\simeq_{(p)}\PU(n)$ for any prime $p$ which does not divide $n$. 
\el

\bq
The quotient map $q$ induces isomorphisms on the $p$-localised homotopy groups (note that $\pi_1(\PU(n))_{(p)}$ is trivial), and hence it is a $p$-local homotopy equivalence by \cite[Theorem 3B(ii)]{hilton-mislin-roitberg}.\qed
\eq

\bl
\label{lem:notdivide}
If the prime $p$ does not divide $n$, then the $p$-primary components of the orders of the Samelson products $\langle \delta_i,1\rangle$ and $\langle \epsilon_i,1\rangle$ coincide.
\el

\bq
Let $p$ be a prime which does not divide $n$. Then $q$ is a $p$-local homotopy equivalence by Lemma \ref{lem:p-equiv}, and hence the commutativity of \eqref{diag} yields
\bse
\langle\delta_i,1\rangle_{(p)}=q^{-1}_{(p)}\circ \langle\epsilon_i,1\rangle_{(p)}\circ(1\land q_{(p)}),
\ese
so the $p$-primary components of the orders of $\langle\delta_i,1\rangle$ and $\langle\epsilon_i,1\rangle$ coincide.\qed
\eq

Hence, when $n$ is prime, the orders of $\langle\delta_i,1\rangle$ and $\langle\epsilon_i,1\rangle$ coincide up to at most their $n$-primary component.

\bl
\label{lem:qiso}
For any $n$, the quotient map $q\cl \SU(n)\to \PU(n)$ induces an isomorphism
\bse
q_*\cl [S^{2i-1}\land \SU(n),\SU(n)]\to[S^{2i-1}\land \SU(n),\PU(n)].
\ese
\el

\bq
Recall that $q\cl \SU(n)\to \PU(n)$ fits into a homotopy fibration sequence
\bse
\cdots\longrightarrow\Z/n\Z\longrightarrow\SU(n)\xrightarrow{\ q\ } \PU(n) \longrightarrow \B(\Z/n\Z).
\ese
Since $\Z/n\Z$ is discrete, applying the functor $[S^{2i-1}\land \SU(n),-]$ yields
\bse
\cdots\longrightarrow 0\longrightarrow[S^{2i-1}\land \SU(n),\SU(n)]\xrightarrow{\ q_*\ }  [S^{2i-1}\land \SU(n),\PU(n)] \longrightarrow 0, 
\ese
whence the statement. \qed
\eq

\bl
\label{cor:atleast}
The order of $\langle\delta_i,1\rangle$ divides the order of $\langle\epsilon_i,1\rangle$. 
\el

\bq
Let $p$ be a prime. If $p^k$ divides the order of $\langle\delta_i,1\rangle$ for some $k\geq 1$, then $p^k$ also divides the order of the composite  $q\circ\langle\delta_i,1\rangle_{(p)}$ by Lemma \ref{lem:qiso}. It then follows, by the commutativity of \eqref{diag}, that the order of $\langle\epsilon_i,1\rangle_{(p)}$ is at least $p^k$. \qed
\eq

For the remainder of this section, we shall restrict to considering $\PU(n)$ when $n$ is an odd prime $p$.

Since $\SU(p)$ is the universal cover of $\PU(p)$ and $H_*(\SU(p);\Z)$ is torsion-free, by \cite[Theorem 1.1]{kishimoto08} we have the following decomposition of $\PU(p)$.

\bl
\label{lem:pup decmop}
For an odd prime $p$, there is a $p$-local homotopy equivalence
\bse
\PU(p)\simeq_{(p)}L\times\prod_{j=2}^{p-1}S^{2j-1}
\ese
where $L$ is an H-space with $\pi_1(L)\cong\Z/p\Z$.\qed
\el

Let $\alpha\cl L_{(p)}\to \PU(p)_{(p)}$ be the inclusion. Then we can write the equivalence of Lemma \ref{lem:pup decmop} as
\bse
L_{(p)}\times\prod_{j=2}^{p-1} S^{2j-1}_{(p)}\xrightarrow{\ \alpha\times \prod_j {\epsilon_j}_{(p)}\ } \bigl(\PU(p)_{(p)}\bigr)^{p-1}\xrightarrow{\ \mu\ }\PU(p)_{(p)},
\ese
where $\mu$ is the group multiplication in $\PU(p)_{(p)}$. We note that this composite is equal to the product
\bse
(\alpha\circ\pr_1)\cdot\prod_{j=2}^{p-1}({\epsilon_j}_{(p)}\circ\pr_j)
\ese
in the group $[L_{(p)}\times\prod_{j=2}^{p-1} S^{2j-1}_{(p)},\PU(p)_{(p)}]$, where $\pr_j$ denotes the projection onto the $j$th factor.

\bl
\label{lem:product}
With the above notation, the localised Samelson product
\bse
\langle \epsilon_i,1\rangle_{(p)} \cl S^{2i-1}_{(p)}\land \PU(p)_{(p)}\to \PU(p)_{(p)}
\ese
is trivial if, and only if, each of $\langle{\epsilon_i}_{(p)},\alpha\rangle$ and $\langle \epsilon_i,\epsilon_j\rangle_{(p)}$, for $2\leq j\leq p-1$, are trivial.
\el

\bq
By \cite[Lemmas 3.3 and 3.4]{hasui16}, $\langle \epsilon_i,1\rangle_{(p)}$ is trivial if, and only if, both $\langle{\epsilon_i}_{(p)},\alpha\rangle$ and $\langle \epsilon_i,\prod_j\epsilon_j\rangle_{(p)}$ are trivial. Applying the same lemmas to the second factor a further $p-3$ times gives the statement. \qed
\eq

We therefore calculate the groups $[S^{2i-1}\land L,\PU(p)]_{(p)}$ and, for $2\leq j\leq p-1$, the homotopy groups $\pi_{2i+2j-2}(\PU(p))_{(p)}$ in order to get an upper bound on the order of $\langle\epsilon_i,1\rangle_{(p)}$.

\bl
\label{lem:pi2i+2j-2}
For $2\leq i\leq p$ and $2\leq j\leq p-1$, the group $\pi_{2i+2j-2}\bigl(\PU(p)\bigr)_{(p)}$ has exponent at most $p$. 
\el

\bq
Decompose $\PU(p)$ as in Lemma \ref{lem:pup decmop}. Observe that, by \cite[Proposition 2.2]{kishimoto08}, we have $\pi_n(L)\cong\pi_n(S^{2p-1})$ for $n\geq 2$,  and hence
\bse
\pi_{2i+2j-2}\bigl(\PU(p)\bigr)_{(p)}\cong\pi_{2i+2j-2}\biggl(L\times\prod_{k=2}^{p-1}S^{2k-1}\biggr)_{\!(p)}\!\!\cong\bigoplus_{k=2}^p\pi_{2i+2j-2}(S^{2k-1})_{(p)}.
\ese
By Toda \cite[Theorem 7.1]{iterated-ii}, if $k\geq 2$ and $r<2p(p-1)-2$, the $p$-primary component of $\pi_{(2k-1)+r}(S^{2k-1})$ is either $0$ or $\Z/p\Z$. Since $2i+2j-2\leq 4p-4$ and
\bse
4p-4<2p(p-1)-2+(2k-1)
\ese
for all $k\geq 2$, the statement follows. \qed
\eq

For the next part of our calculation, we will need a certain mod-$p$ decomposition of $\Sigma^2 L$ which will, in turn, require some cohomological information. The mod-$p$ cohomology algebra of $\PU(n)$, with $p$ any prime and $n$ arbitrary, was determined by Baum and Browder in \cite[Corollary 4.2]{baum-browder}. In particular, we have:
\bl
\label{lem:pup-cohom}
For $p$ an odd prime, there is an algebra isomorphism
\begin{equation*}
    H^*(\PU(p);\Z/p\Z)\cong \Lambda(x_1,x_3,\ldots,x_{2p-3})\otimes \frac{\Z/p\Z[y]}{(y^p)},
\end{equation*}
with $|x_d|=d$, $|y|=2$, and $\beta(x_1)=y$.\qed 
\el

For $m\geq 2$, denote by $P^{m}(p)$ the mod-$p$ Moore space defined as the homotopy cofibre of the degree $p$ map 
\bse
S^{m-1}\xrightarrow{\ p \ } S^{m-1}\longrightarrow P^{m}(p)
\ese
on the sphere $S^{m-1}$. In other words, $P^m(p)=S^{m-1}\cup_{p} e^m$. Note that, by extending the cofibre sequence to the right, we see that $\Sigma P^{m}(p)\simeq P^{m+1}(p)$. 

\bl
\label{lem:susp-L}
For $p$ an odd prime, there is a $p$-local homotopy equivalence
\bse
\Sigma L \simeq_{(p)} A\lor \bigvee_{k=2}^{p-1}P^{2k+1}(p),
\ese
with $H_*(A;\Z/p\Z)$ generated by $\{u,v,w\}$, with $|u|=2$, $|v|=3$, $|w|=2p$, and subject to the relation $\beta(v)=u$.
\el

\bq
By decomposing $\PU(p)$ as in Lemma \ref{lem:pup decmop}, taking mod-$p$ cohomology and comparing with Lemma \ref{lem:pup-cohom}, we obtain
\begin{equation*}
    H^*(L;\Z/p\Z)\cong \Lambda(x_1)\otimes \frac{\Z/p\Z[y]}{(y^p)}.
\end{equation*}
Since $H^*(L;\Z/p\Z)$ is of finite type and self-dual, $H^*(L;\Z/p\Z)\cong H_*(L;\Z/p\Z)$ as Hopf algebras. Then, as $H_*(L;\Z/p\Z)$ is primitively generated and $L$ is a connected H-space (being a retract of $\PU(p)$), by \cite[Theorem 4.1]{cohen1976} there is a decomposition
\bse
\Sigma L \simeq_{(p)} A_1\lor A_2\lor \cdots\lor A_{p-1},
\ese
with each summand $A_j$ having homology $H_*(A_j;\Z/p\Z)$ generated by the suspensions of monomials in $H_*(L;\Z/p\Z)$ of length $j$ (modulo $p-1$), where by length of a monomial one means the number of (not necessarily distinct) factors in that monomial.

Let $\overline{x}_1$ and $\overline{y}$ denote the duals of $x_1$ and $y$, and let $\sigma$ denote the suspension isomorphism for homology.
Then, $H_*(A_1;\Z/p\Z)$ is generated by $\sigma(\overline{x}_1)$, $\sigma(\overline{y})$, and $\sigma(\overline{x}_1\overline{y}^{p-1})$, in degrees 2, 3, and $2p$, respectively. Furthermore, by the stability of the Bockstein operator $\beta$, we also have $\beta(\sigma(\overline{y}))=\sigma(\overline{x}_1)$.

On the other hand, for $j\neq 1$, the homology $H_*(A_j;\Z/p\Z)$ is generated by the elements $\sigma(\overline{x}_1\overline{y}^{j-1})$ and $\sigma(\overline{y}^j)$, in degrees $2j$ and $2j+1$, respectively, subject to the relation $\beta(\sigma(\overline{y}^j))=\sigma(\overline{x}_1\overline{y}^{j-1})$. As the homotopy type of Moore spaces is uniquely characterised by their homology, we must have $A_j\simeq P^{2j+1}(p)$ for $j\neq 1$, yielding the decomposition in the statement. \qed
\eq

With $A$ as in Lemma \ref{lem:susp-L}, we have:

\bl
\label{lem:a-decomp}
There is a $p$-local homotopy equivalence
\bse
\Sigma A \simeq_{(p)} P^4(p)\lor S^{2p+1}.
\ese
\el

\bq
Localise at $p$ throughout. By looking at the degrees of the generators of $H_*(A;\Z/p\Z)$ in Lemma \ref{lem:susp-L}, we see that the 3-skeleton of $A$ is $P^3(p)$.

Let $f\cl S^{2p-1}\to P^3(p)$ be the attaching map of the top cell of $A$, and let $F$ be the homotopy fibre of $\rho\cl P^3(p)\to S^3$, the pinch map to the top cell of $P^3(p)$. As $\pi_{2p-1}(S^3)\cong 0$, the map $f$ lifts through $F$ via some map $\lambda\cl S^{2p-1}\to F$, as in the diagram:
\begin{equation*}
    \begin{tikzcd}
    & F \ar[d]\\
    S^{2p-1}\ar[ur,"\lambda"] \ar[r,"f"]\ar[dr,"*"']& P^3(p)\ar[d,"\rho"]\\
    & S^3.
    \end{tikzcd}
\end{equation*}
Let $j\cl S^3\to P^4(p)$ be the inclusion of the bottom cell and let $S^3\{p\}$ be the homotopy fibre of the degree $p$ map on $S^3$. As $j$ has order $p$, there is a homotopy fibration diagram
\begin{equation*}
    \begin{tikzcd}
    \Omega S^3\ar[r]\ar[d,"\Omega j"] & S^3\{p\}\ar[d,"s"]\ar[r] & S^3 \ar[d]\ar[r,"p"]& S^3\ar[d]\ar[d,"j"]\\
    \Omega P^4(p) \ar[r,equal]&    \Omega P^4(p) \ar[r]& * \ar[r]& P^4(p) 
    \end{tikzcd}
\end{equation*}
which defines a map $s\cl S^3\{p\}\to\Omega P^4(p)$. For connectivity reasons, the suspension map $P^3(p)\xrightarrow{\, E \,}\Omega P^4(p)$ factors as the composite $P^3(p)\xrightarrow{\, \iota\,} S^3\{p\}\xrightarrow{\, s\,} \Omega P^4(p)$, where $\iota$ is the inclusion of the bottom Moore space. Furthermore, there is a homotopy fibration diagram
\begin{equation*}
    \begin{tikzcd}
   F \ar[d] \ar[r] & P^3(p) \ar[d,"\iota"]\ar[r,"\rho"] & S^3\ar[d,equal]\\
   \Omega S^3\ar[r] & S^3\{p\}\ar[r]& S^3.
    \end{tikzcd}
\end{equation*}
Putting this together gives a commutative diagram
\begin{equation*}
    \begin{tikzcd}
    & F\ar[r]\ar[d]&\Omega S^3 \ar[d]\ar[r,equal] &\Omega S^3\ar[d,"\Omega j"] \\
    S^{2p-1}\ar[ur,"\lambda"]\ar[r,"f"] & P^3(3)\ar[r,"\iota"]&S^3\{p\}\ar[r,"s"]&\Omega P^4(p) .
    \end{tikzcd}
\end{equation*}
Thus $E\circ f$ factors through $\Omega j$, implying that $\Sigma f$ factors as the composite
\bse
S^{2p}\xrightarrow{\ \hat f\ } S^3 \xrightarrow{\ j \ } P^4(p)
\ese
for some map $\hat f$.

As $\pi_{2p}(S^3) \cong\Z/p\Z$, we must have $\hat f=t\alpha$, where $\alpha$ is a generator of $\pi_{2p}(S^3)$. If $\Sigma f$ were essential, then $t\neq 0$. However, the element $\alpha$ would then be detected by the Steenrod operation $\cP^1$ in the cohomology of $\Sigma A$. This would, in turn, imply that $\cP^1$ were non-trivial in $H^*(A;\Z/p\Z)$,  and hence in $H^*(\PU(p);\Z/p\Z)$. However, $\cP^1(H^*(\PU(p);\Z/p\Z))=0$, and thus we must have had $\Sigma f \simeq *$. \qed 
\eq

\bl
\label{lem:suspL}
The exponent of the group $[S^{2i-1}\land L,\PU(p)]_{(p)}$ is at most $p$. 
\el

\bq
By the decompositions in Lemmas \ref{lem:susp-L} and \ref{lem:a-decomp}, we have
\begin{align*}
    [S^{2i-1}\land L,\PU(p)]_{(p)}  &\cong \bigl[S^{2i-2}\land  \bigl(A\lor \textstyle\bigvee_{k=2}^{p-1}P^{2k+1}(p)\bigr),\PU(p)\bigr]_{(p)}\\
    &\cong \bigl[S^{2i-3}\land  \bigl(S^{2p+1}\lor \textstyle\bigvee_{k=1}^{p-1}P^{2k+2}(p)\bigr),\PU(p)\bigr]_{(p)}\\
    & \cong \pi_{2i+2p-2}(\PU(p))_{(p)}\oplus\bigoplus_{k=1}^{p-1}[P^{2k+2i-1}(p),\PU(p)]_{(p)}.
\end{align*}
Since $2i+2p-2\leq 4p-2<2p(p-1)+1$ for $p\geq 3$, the group $\pi_{2i+2p-2}(\PU(p))_{(p)}$ consists of elements of order at most $p$ by the same argument as in Lemma \ref{lem:pi2i+2j-2}.

On the other hand, by \cite[Theorem 7.1]{neisendorfer80}, the groups $[P^{2k+2i-1}(p),\PU(p)]$ have exponent at most $p$ (since, for $m\geq 3$, the identity on $P^m(p)$ has order $p$), whence the statement. \qed
\eq

Lemmas \ref{lem:product}, \ref{lem:pi2i+2j-2}, and \ref{lem:suspL} together imply:

\bl
\label{cor:atmost}
The order of the Samelson product
\bse
\langle \epsilon_i,1\rangle_{(p)}\cl S^{2i-1}_{(p)}\land\PU(p)_{(p)}\to \PU(p)_{(p)}
\ese
is at most $p$. \qed
\el

We now have all the ingredients necessary to prove Theorem \ref{thm:pup-samels}.

\bq[of Theorem \ref{thm:pup-samels}]
Consider the following commutative diagram
\bse
\begin{tikzcd}[row sep=large,column sep=large]
S^{2i-1}\land S^{2(p-i)+1} \ar[r,"{\langle\eta_i,\eta_{p-i-1}\rangle}"] \ar[d,"1\land \delta_{p-i-1}"]&[12pt] \U(p)\\
S^{2i-1}\land \SU(p) \ar[r,"{\langle\delta_i,1\rangle}"] & \SU(p)\ar[u,"\iota"']
\end{tikzcd}
\ese
where $\iota\cl \SU(p)\to\U(p)$ is the inclusion and $\eta_i:=\iota_*(\delta_i)$. By the unnumbered corollary in Bott \cite[p. 250]{bott60}, the map $\langle\eta_i,\eta_{p-i-1}\rangle$ is non-trivial and $p$ divides its order. Hence, the order of $\langle\delta_i,1\rangle_{(p)}$ is at least $p$. The result now follows from Lemmas \ref{lem:notdivide}, \ref{cor:atleast} and \ref{cor:atmost}. \qed
\eq

\section{Homotopy invariants of $\PU(n)$-gauge groups}
\label{sec:hinv}

The content of Lemma \ref{thm:sutopu} is a straightforward observation about how certain homotopy invariants of $\SU(n)$-gauge groups relate to the corresponding invariants of $\PU(n)$-gauge groups.

\bl
\label{thm:sutopu}
Let $n$ be arbitrary and $X$ be a simply-connected space. Suppose further that we have $[X,\SU(n)]\cong 0$. Then, the quotient map $q\cl \SU(n)\to\PU(n)$ induces an isomorphism of groups
\begin{equation*}
    [X,\cG_{i,k}(\SU(n))]\cong [ X,\cG_{i,k}(\PU(n))]
\end{equation*}
for any $2\leq i\leq 2n$ and any $k\in \Z$.
\el

\bq
Since $[ X,\SU(n)]\cong 0$ and $X$ is simply-connected, applying the functor $[X,-]$ to the homotopy fibration sequence
\bse
\cdots\longrightarrow\Z/n\Z\longrightarrow\SU(n)\xrightarrow{\ q\ } \PU(n) \longrightarrow \B(\Z/n\Z)
\ese
shows that $[X,\PU(n)]\cong 0$ also. 

Applying now the functor $[\Sigma X,-]$ to the homotopy fibration sequence
\bse
\PU(n)\xrightarrow{\ \p_{i,k}\ }\Omega_0^{2i-1}\PU(n)\longrightarrow \B\cG_{i,k}(\PU(n))\longrightarrow\B\PU(n)
\ese
described in Section \ref{sec:gauge}, as well as to its $\SU(n)$ analogue, yields the following commutative diagram
\bse
\begin{tikzcd}[column sep=normal]
{[\Sigma X,\SU(n)]} \ar[r,"(\p_{i,k})_*"]\ar[d,"q_*"] &[9pt] {[\Sigma^{2i}X,\SU(n)]} \ar[r]\ar[d,"q_*"] & {[ X,\cG_{i,k}(\SU(n))]} \ar[r]\ar[d] & \ar[d] 0\\
{[\Sigma X,\PU(n)]} \ar[r,"(\p_{i,k})_*"] & {[\Sigma^{2i}X,\PU(n)]} \ar[r] & {[ X,\cG_{i,k}(\PU(n))]} \ar[r] & 0
\end{tikzcd}
\ese
where the rows are exact and the two leftmost vertical maps are isomorphisms by the same argument as in the proof of Lemma \ref{lem:qiso}. The statement now follows from the five lemma. \qed
\eq

Hamanaka and Kono showed in \cite[Theorem 1.2]{hamanaka-kono} that, for principal $\SU(n)$-bundles over $S^4$, the homotopy equivalence $\cG_{2,k}(\SU(n))\simeq \cG_{2,l}(\SU(n))$ implies that $(n(n^2-1),k)=(n(n^2-1),l)$. As an application of Lemma \ref{thm:sutopu}, let us show that the analogue of this result holds for $\PU(n)$-gauge groups.

\bc
\label{cor:sutopus4}
Let $n> 3$. For principal $\PU(n)$-bundles over $S^4$,if
\bse
\cG_{2,k}(\PU(n))\simeq \cG_{2,l}(\PU(n)),
\ese
then $(n(n^2-1),k)=(n(n^2-1),l)$.
\ec

\bq
First, suppose that $n$ is even. Note that we have
\bse
\pi_{2n-4}(\SU(n))\cong\pi_{2n-2}(\SU(n))\cong 0.
\ese
Hence, applying Lemma \ref{thm:sutopu} with $X=S^{2n-4}$ and $X=S^{2n-2}$, we find
\bse
\pi_{2n-4}\bigl(\cG_{2,k}(\PU(n))\bigr)\cong\pi_{2n-4}\bigl(\cG_{2,k}(\SU(n))\bigr)
\ese
and
\bse
\pi_{2n-2}\bigl(\cG_{2,k}(\PU(n))\bigr)\cong\pi_{2n-2}\bigl(\cG_{2,k}(\SU(n))\bigr).
\ese
So the result follows for $n$ even by \cite[Proposition 4.2]{sutherland92}.

When $n$ is odd, we have from \cite{hamanaka-kono} that $[\Sigma^{2n-6} \C P^2,\SU(n)]\cong 0$. Hence, applying Lemma \ref{thm:sutopu} with $X=\Sigma^{2n-6}\C P^2$, we find
\bse
[\Sigma^{2n-6} \C P^2,\cG_{2,k}(\PU(n))]\cong [\Sigma^{2n-6} \C P^2,\cG_{2,k}(\SU(n))].
\ese
So the result follows for $n$ odd by \cite[Corollary 2.6]{hamanaka-kono}.\qed
\eq

Following the work of \cite{hamanaka-s6}, Mohammadi and Asadi-Golmankhaneh \cite{su(n)-s6} recently showed that, for $\SU(n)$-bundles over $S^6$, an equivalence $\cG_{3,k}(\SU(n))\simeq\cG_{3,l}(\SU(n))$ implies that 
\bse
\bigl((n - 1)n(n + 1)(n + 2), k\bigr) = \bigl((n - 1)n(n + 1)(n + 2), l\bigr).
\ese
Hence, we also have:

\bc
\label{cor:s6}
Let $n\geq 3$. For principal $\PU(n)$-bundles over $S^6$, if there is a homotopy equivalence $\cG_{3,k}(\PU(n))\simeq\cG_{3,l}(\PU(n))$, then
\bse
((n - 1)n(n + 1)(n + 2), k) = ((n - 1)n(n + 1)(n + 2), l).
\ese
\ec

\bq
Apply Lemma \ref{thm:sutopu} with $X=\Sigma^{2n-6}\C P^2$ and the result of \cite{su(n)-s6}. \qed


\eq

\section{Special cases}

\subsection{$\PU(p)$-bundles over $S^4$} Theriault showed in \cite{theriault-sun} that, after localisation at an odd prime $p$ and provided $n<(p-1)^2+1$, the order of the Samelson product $\langle\delta_2,1\rangle\cl S^3\land \SU(n)\to \SU(n)$ is the $p$-primary component of the integer $n(n^2-1)$. 
It then follows immediately from Theorem \ref{thm:pup-samels} that:

\bc
After localisation at an odd prime, the order of the Samelson product $\langle\epsilon_2,1\rangle\cl S^3\land \PU(p)\to \PU(p)$ is $p(p^2-1)$. \qed
\ec

\subsection{$\PU(5)$-bundles over $S^4$} In \cite{theriault-su5}, Theriault showed that the order of $\langle \delta_2,1\rangle\cl S^3\wedge \SU(5)\to \SU(5)$ is 120. Hence, by Theorem \ref{thm:pup-samels}, the order of $\langle \epsilon_2,1\rangle\cl S^3\wedge \PU(5)\to \PU(5)$ is also 120.

\bq[of Theorem \ref{thm:pu5}]
Part (i) follows from Corollary \ref{cor:sutopus4}, while part (ii) follows from Lemma \ref{lem:theriault}. \qed
\eq

\subsection{$\PU(3)$-bundles over $S^6$} Hamanaka and Kono showed in \cite{hamanaka-s6} that the order of $\langle\delta_3,1\rangle\cl S^5\land \SU(3)\to \SU(3)$ is $120$. It follows from Theorem \ref{thm:pup-samels} that the order of $\langle\epsilon_3,1\rangle\cl S^5\land \PU(3)\to\PU(3)$ is also $120$.

\bq[of Theorem \ref{thm:pu3}]
As the homotopy groups $\pi_n(\Omega_0^5\PU(3))\cong\pi_{n+5}(\PU(3))$ are all finite, the ``if'' direction follows from Lemma \ref{lem:hamanaka}, while the ``only if'' direction follows from Corollary \ref{cor:s6}. \qed
\eq


%
%

\bibliographystyle{spmpsci}      

\bibliography{references}


\end{document}

%% file: sc.tex
  \def\bc{\begin{corollary}}
  \def\ec{\end{corollary}}
  \def\bd{\begin{definition}}
  \def\ed{\end{definition}}
  \def\bse{\begin{equation*}}
  \def\ese{\end{equation*}}
  \def\be{\begin{example}}
  \def\ee{\end{example}}
  \def\bi{\begin{IEEEeqnarray*}}
  \def\ei{\end{IEEEeqnarray*}}
  \def\bl{\begin{lemma}}
  \def\el{\end{lemma}}
  \def\bp{\begin{proposition}}
  \def\ep{\end{proposition}}
  \def\bq{\begin{proof}}
  \def\eq{\end{proof}}
  \def\bt{\begin{theorem}}
  \def\et{\end{theorem}}

  \def\B{\mathrm{B}}
  \def\C{\mathbb{C}}
  \def\cG{\mathcal{G}}
  \def\cP{\mathcal{P}}
  \def\cl{\colon}
  \def\Map{\mathrm{Map}}

  \def\PSU{\mathrm{PSU}}
  \def\PU{\mathrm{PU}}
  \def\p{\partial}
  \def\pr{\mathrm{pr}}

  \def\SU{\mathrm{SU}}

  \def\U{\mathrm{U}}

  \def\Z{\mathbb{Z}}

\DeclareMathOperator{\ev}{ev}

\DeclareMathOperator{\SO}{SO}